\documentclass[12pt]{article}

\usepackage{diagbox}
\usepackage{mathtools}
\usepackage{bbm}
\usepackage{latexsym}
\usepackage{epsfig}
\usepackage{amsmath,amsthm,amssymb,enumerate,hyperref}
\usepackage[active]{srcltx}
\usepackage{color}

\newcounter{rot}

\newcommand{\hatt}[1]{\widehat #1}

\def\hS{\hatt{S}}

\def\a{\alpha} \def\b{\beta} \def\d{\delta} 
    \def\g{\gamma}

 \def\m{\mu}  
   
 \def\om{\omega}

\def\hY{\hatt{Y}}

\newtheorem{theorem}{Theorem}
\newtheorem{lemma}[theorem]{Lemma}

\newcommand{\proofstart}{{\bf Proof\hspace{2em}}}
\newcommand{\proofend}{\hspace*{\fill}\mbox{$\Box$}}

\newcommand{\brac}[1]{\left(#1\right)}

\newcommand{\bfrac}[2]{\left(\frac{#1}{#2}\right)}

\newcommand{\set}[1]{\left\{#1\right\}}

\def\Pr{mathbb{P}}

\newcommand{\ignore}[1]{}

\newcommand{\beq}[2]{\begin{equation}\label{#1}#2\end{equation}}

\newcommand{\multstar}[1]{\begin{multline*}#1\end{multline*}}

\usepackage{tikz}
\usetikzlibrary{arrows}
\usetikzlibrary{decorations}
\usetikzlibrary{shapes.misc}

\def\am{\a_{\min}}

\def\hS{\widehat{X}}
\def\hT{\widehat{Y}}

\begin{document}
\author{Alan Frieze\\\ Department of Mathematical Sciences,\\ Carnegie Mellon University,\\ Pittsburgh PA15213,\\ USA.\\alan@random.math.cmu.edu.\thanks{Research supported in part by NSF grant DMS1661063}}

\title{A note on randomly colored matchings in random bipartite graphs}
\maketitle

\begin{abstract}{}
We are given a bipartite graph that contains at least one perfect matching and where each edge is colored from a set $Q=\set{c_1,c_2,\ldots,c_q}$. Let \\
$Q_i=\set{e\in E(G):c(e)=c_i}$, where $c(e)$ denotes the color of $e$. The perfect matching color profile $mcp(G)$ is defined to be the set of vectors $(m_1,m_2,\ldots,m_q)\in [n]^q$ such that there exists a perfect matching $M$ such that $|M\cap Q_i|=m_i$. We give bounds on the matching color profile for a randomly colored random bipartite graph.
\end{abstract}
\section{Introduction}
We consider the following problem: we are given a random bipartite graph $G$ in which each edge is given a random color from a set $Q=\set{c_1,c_2,\ldots,c_q}$. An edge $e$ is colored $c(e)=c_i$ with probability $\a_i$ where $\a_i>0$ is a constant.  Let $Q_i=\set{e\in E(G):c(e)=c_i}$, where $c(e)$ denotes the color of $e$. The perfect matching color profile $mcp(G)$ is defined to be the set of vectors $(m_1,m_2,\ldots,m_q)\in [n]^q$ such that there exists a perfect matching $M$ such that $|M\cap Q_i|=m_i$. We give bounds on the matching color profile for a randomly colored random bipartite graph.

Randomly colored random graphs have been studied recently in the context of (i) rainbow matchings and Hamilton cycles, see for example \cite{BF}, \cite{CF}, \cite{FK}, \cite{JW}; (ii) rainbow connection see for example \cite{DFT15}, \cite{gnprainbow}, \cite{HR}, \cite{M}, \cite{KKS}; (iii) pattern colored Hamilton cycles, see for example \cite{AF}, \cite{EFK}. This paper can be considered to be a contribution in the same genre. One can imagine a possible interest in the color profile via the following scenario: suppose that $A$ is a set of tools and $B$ is a set of jobs where edge $\set{a,b}$ indicates that $b$ can be completed using $a$. If colors represent people, then one might be interested in equitably distributing jobs. I.e. determining whether $(n/q,n/q,\ldots,n/q)\in mcp(G)$. In any case, we find the problem interesting.

We will consider $G$ to be the random bipartite graph $G_{n,n,p}$ where $p=\frac{\log n+\om}{n},\,\om=\om(n)\to\infty$ where $\om=o(\log n)$. Erd\H{os} and R\'enyi \cite{ER} proved that $G$ has a perfect matching w.h.p. We will prove the following theorem: let $\a_1,\a_2,\ldots,\a_q,\b$ be positive constants such that $\a_1+\a_2+\cdots+\a_q=1$ and $\b<1/q$. Let
\[
\am=\min\set{\a_i:i\in [q]}.
\]
\begin{theorem}\label{th1}
Let $G$ be the random bipartite graph $G_{n,n,p}$ where $p=\frac{\log n+\om}{n},\,\om=\om(n)\to\infty$ where $\om=o(\log n)$. Suppose that the edges of $G$ are independently colored with colors from $C=\set{c_1,c_2,\ldots,c_q}$ where $\Pr(c(e)=c_i)=\a_i$ for $e\in E(G),i\in[q]$. Let $m_1,m_2,\ldots,m_q$ satisfy: (i) $m_1+\cdots+m_q=n$ and (ii) $m_i\geq \b n,i\in [q]$. Then w.h.p., there exists a perfect matching $M$ in which exactly $m_i$ edges are colored with $c_i,i=1,2,\dots,q$.
\end{theorem} 
It is clear that w.h.p. $(n,0,\ldots,0)\notin mcp(G)$. This is because the bipartite graph induced by edges of color $c_1$ is distributed as $G_{n,n,\a_1p}$ and this contains isolated vertices w.h.p. On the other hand, if $p\geq \frac{q(\log n+\om)}{\am n}$ then w.h.p. $mcp(G)=[n]^q$. To see this, suppose that $m_1\leq m_2\leq \cdots m_q\leq n$. Suppose we have found a matching that uses $m_i$ edges of color $c_i$ for $i\geq 0$. Let $n'=n-m_1-\cdots-m_i$. Then the random bipartite graph induced by vertices not in $M$ and having edges of color $c_i$ has density at least $\frac{q\a_in'}{\am n}\cdot\frac{\log n+\om}{n'}\geq \frac{\log n'+\om/2}{n'}$ and so has a perfect matching w.h.p.

{\bf Open Question:} What is the threshold for $mcp(G)=[0,n]^q$?
\section{Structural Lemma}
Suppose that the bipartition of $V(G)$ is denoted $A,B$. For sets $S\subseteq A,T\subseteq B$ we let $e_i(S,T)$ denote the number of $S:T$ edges of color $c_i$. We say that vertex $u$ is $c_i$-adjacent to vertex $v$ if the edge $\set{u,v}$ exists and has color $c_i$.
\begin{lemma}\label{lem1}
 Let $p=\frac{\log n+\om}{n},\,\om=\om(n)\to\infty$ where $\om=o(\log n)$. Then w.h.p.
\begin{enumerate}[(a)]
\item $S\subseteq A, T\subseteq B$ and $\g_a\log n\leq |S|\leq n_0=\g_a n/\log n$ and $|T|\leq\a_i \eta|S|\log n$ where $\g_a=\eta/(20\a_i)$ implies that $e_i(S:T)\leq 2\a_i\eta  |S|\log n$ for $i=1,2,\ldots,q$.
\item There do not exist sets $X\subseteq S\subseteq A, T\subseteq B$ and $i\in [q]$ such that $|S|,|T|\geq \beta n$ and $|X|=\g_b|S|/\log n,\g_b=10\log(e/\b)/\a_i$ and such that each $x\in X$ is $c_i$-adjacent to fewer than $\a_i\b\log n/10$ vertices in $T$.
\item  There do not exist sets $X\subseteq S\subseteq A, T\subseteq B$ and $i\in [q]$ such that $|S|,|T|\geq \beta n$ and $|X|=|S|/\log n$ and a set $Z\subseteq T,|Z|=\g_bn/\log n$ such that each $x\in X$ is $c_i$-adjacent to $k=\frac{10\log n}{\log\log n}$ vetices in $Z$.
\item There do not exist sets $S\subseteq A, T\subseteq B$ and $i\in [q]$ such that $|S|,|T|\geq \beta n$ such that there are more than $\g_dn/\log n,\g_d=\frac{4}{\a_i}\log\bfrac{e}{\b}$ vertices in $T$ that not $c_i$-adjacent to a vertex in $S$.
\item Fix $\g,\d>0$ constants. Then w.h.p. there do not exist sets $S,T$ with $|S|=|T|=\g n/\log n$ such that $e_i(S,T)\geq \d|S|\log n/\log\log n$.
\item There do not exist sets $S\subseteq A, T\subseteq B$ and $i\in [q]$ such that $|S|,|T|\geq \beta n/10$ such that $e_i(S,T)=0$.
\end{enumerate}
\end{lemma}
\proofstart\\
(a) The probability that the condition is violated can be bounded by
\multstar{
\sum_{s=\g_a\log n}^{n_0}\sum_{t=1}^{\a_i\eta s\log n}\binom{n}{s}\binom{n}{t}\binom{st}{2\a_i\eta s\log n}(\a_ip)^{2\a_i\eta s\log n}\\
\leq \sum_{s=\g_a\log n}^{n_0}\sum_{t=1}^{\a_i\eta s\log n}\bfrac{ne}{s}^s\bfrac{ne}{t}^t\bfrac{est\a_ip}{2\a_i\eta s\log n}^{2\a_i\eta s\log n}\\
\leq \sum_{s=\g_a\log n}^{n_0}\sum_{t=1}^{\a_i\eta s\log n}\bfrac{ne}{s}^s\bfrac{ne}{\a_i\eta s\log n}^{\a_i\eta s\log n}\bfrac{etp}{2\eta \log n}^{2\a_i\eta s\log n}\\
\leq\sum_{s=\g_a\log n}^{n_0}\sum_{t=1}^{\a_i\eta s\log n}\brac{\bfrac{ne}{s}^{1/2\a_i\eta\log n} \bfrac{ne}{\a_i\eta s\log n}^{1/2}\cdot \frac{e^{1+o(1)}\a_is\log n}{2n}}^{2\a_i\eta s\log n}\\
\leq \sum_{s=\g_a\log n}^{n_0}\sum_{t=1}^{\a_i\eta s\log n}\bfrac{s\log n}{n}^{\a_i\eta s\log n-s}(\log n)^{s}\bfrac{e^{3/2+o(1)}\a_i^{1/2}}{2\eta^{1/2}}^{2\a_i\eta s\log n}=o(1).
}
(b) The probability that the condition is violated can be bounded by
\[
\binom{n}{\b n}^2\binom{\b n}{\g_bn/\log n}\brac{e^{-\a_i\b/4}}^{\g_bn}\leq \brac{\bfrac{e}{\b}^{(2\b+o(1))n}e^{-\a_i\b\g_b/4}}^n=o(1).
\]
The factor $e^{-\a_i\b\g_b/4}$ comes from applying a Chernoff bound.

(c) We can assume w.l.o.g. that $|S|=|T|=\b n$. The probability that the condition is violated can be bounded by
\multstar{
\binom{n}{\b n}^2\binom{\b n}{n/\log n}\binom{\b n}{\g_bn/\log n}\brac{\binom{\g_bn/\log n}{k}(\a_ip)^k}^{n/\log n}\\
\leq \bfrac{e}{\b}^{(2\b+o(1))n}\bfrac{e\g_b\a_i}{k}^{kn/\log n}=o(1).
}
(d) The probability that the condition is violated can be bounded by
\[
\binom{n}{\b n}^2\binom{\b n}{\g_dn/\log n}(1-\a_ip)^{\b n\g_dn/\log n}\leq \brac{\bfrac{e}{\b}^{2+o(1)} e^{-\a_i\g_d}}^{\b n}=o(1).
\]
(e) The probability that the condition is violated can be bounded by
\multstar{
\binom{n}{\g n/\log n}^2\binom{\g^2n^2/(\log n)^2}{\d n/\log\log n}p^{\d n/\log\log n} \leq\\ \bfrac{e\log n}{\g}^{2\g n/\log n} \bfrac{\g^2e\log\log n}{\d \log n}^{\d n/\log\log n}=o(1).
}
(f) The probability that the condition is violated can be bounded by
\[
2^{2n}(1-p)^{\b^2n^2/100}=o(1).
\]
\proofend
\section{Proof of Theorem \ref{th1}}
\proofstart
Assume from now on that the high probability conditions of Lemma \ref{lem1} are in force. Let $M$ be a perfect matching and let $\m_i=|M\cap Q_i|$ for $i\in [q]$. Suppose that $\m_1>m_1\geq \b n$ and $\b n\leq \m_2<m_2$. We show that we can find another matching $M'$ such that $|M'\cap Q_1|=\m_1-1$ and $|M'\cap Q_2|=\m_2+1$. We do this by finding an alternating cycle with edge sequence  $C=(e_1,f_1,\ldots,e_\ell,f_\ell)$ and vertex sequence $(x_1\in A,y_1\in B,x_2,\ldots,x_\ell,y_\ell,x_1)$ such that (i) $e_i=\set{x_i,y_i}\in M,$, (ii) $f_i=\set{y_i,x_{i+1}}\notin M, i\in[\ell]$, (iii) $e_1\in Q_1$ and (iv) $E(C)\setminus \set{e_1}\subseteq Q_2$. Repeating this for pairs of colors, one over-subscribed and one under-subscribed we eventually achieve our goal.
It is sufficient to consider this case, seeing as we can always w.h.p. find a matching that has been randomly colored with $\approx \a_in$ edges of color $c_i$, $i=1,2,\ldots,q$.

Next let $A_i=V(M\cap Q_i)\cap A$ and $B_i=V(M\cap Q_i)\cap B$ for $i\in [q]$ and for $S\subseteq A$ let $N_i(S)=\set{b\in B:\exists a\in S\ s.t.\ . \set{a,b}\in Q_i}$ and $N_i(a)=N_i(\set{a})$. Then let
\begin{align*}
D_0'&=\set{a\in A_2:|N_2(a)\cap B_2|\geq \frac{\a_2\b\log n}{10}}.\\
D_0&=\set{a\in A_1:|N_2(a)\cap M(A_2\setminus D_0')|\leq k_0=\frac{10\log n}{\log\log n}}.
\end{align*}
It follows from Lemma \ref{lem1}(b) that 
\[
|M(A_2\setminus D_0')|\leq \frac{\g_bn}{\log n}.
\]
It then follows from Lemma \ref{lem1}(c) that if $W_0=A_1\setminus D_0$ then
\beq{W0}{
|W_0|\leq \frac{n}{\log n}.
}
We now define a sequence of sets $W_0,W_1,\ldots$ where $W_{j+1}$ is obtained from $W_j$ by adding a vertex of $A_2\setminus W_j$ for which $|N_2(a)\cap M(W_j)|\geq k_0$. Now consider $S=W_t,T=M(W_t)$ for some $t\geq 1$. Then we have
\[
|S|=|T|\leq t+\frac{n}{\log n}\text{ and }e_2(S,T)\geq tk_0.
\]
Given Lemma \ref{lem1}(e) with $\d=5,\g=2$, we see that this sequence stops with $t=t^*\leq 4n/\log n$. So we now let $R_0=A_2\setminus W_{t^*}$. We note that
\beq{props}{
\begin{split}
|R_0|\geq \b n-\frac{5n}{\log n}\\
a\in R_0\text{ implies }|N_2(a)\cap M(R_0)|\geq \frac{\a_2\b\log n}{10}-k_0.
\end{split}
}
We now fix some $a_0\in R_0$ and define a sequence of sets $X_0,Y_0,X_1,Y_1,\ldots$ where $X_j\subseteq R_0$ and $Y_j\subseteq B_2$. We let $X_0=\set{a_0}$ and then having defined $X_i,i\geq 0$ we let
\[
Y_{i}=N_2(X_i)\text{ and }X_{i+1}=\brac{M^{-1}(Y_i)\setminus \bigcup_{j\leq i}X_j}\cap R_0.
\]
We claim that for  $i\geq 0$,
\beq{cl1}{
|X_i|\leq \frac{n}{200\log n}\text{ implies that }|X_{i+1}|\geq \frac{\a_2\b\log n}{25}|X_i|.
}
We verify \eqref{cl1} below. Assuming its truth, there exists a smallest $k$ such that
\beq{Sk}{
|X_k|\geq \frac{\a_2\b n}{5000}.
}

Starting with $\hT_0=\set{b_0}$ where $b_0=M(a_0)\in \hat{R}_0$, we can similarly construct a sequence of sets $\hT_1,\hS_1,\ldots$  where $\hS_j\subseteq M^{-1}(\hat{R}_0)$ and $\hY_j\subseteq\hat{R}_0$. Here $\hat{R}_0$ is the equivalently defined set to $R_0$ in $B_2$. We can assume that $b_0\in \hat{R}_0$, because of the sizes of the sets $R_0,\hat{R}_0$. More precisely, by \eqref{W0}, there will be $o(n)$ choices for $a_0$ for which $b_0\notin \hat{R}_0$. Having defined $\hT_i$ we let
\[
\hS_{i}=N_2(\hT_i)\text{ and }\hY_{i+1}=\brac{M(\hS_i)\setminus \bigcup_{j\leq i}\hT_j}\cap \hat{R}_0.
\]
and then let $\hT_{i+1}=M(\hS_i)$. The equivalent of \eqref{cl1} will be 
\beq{cl2}{
|\hT_i|\leq \frac{n}{200\log n}\text{ implies that }|\hT_{i+1}|\geq \frac{\a_2\b\log n}{25}|\hT_i|.
}
Assuming its truth, there exists $\ell$ such that 
\beq{Tell}{
|\hT_\ell|\geq \frac{\a_2\b n}{5000}.
}
It follows from Lemma \ref{lem1}(f) that at least 9/10 of the vertices of $A_1$ have a $c_2$-neighbor in $\hat{T}_0$ and at least 9/10 of the vertices of $B_1$ have a $c_2$-neighbor in $R_0$. We deduce from this that there is a pair $x_0\in A_1,y_0=M(x_0)\in B_1$ such that $N_2(x_0)\cap \hT_\ell\neq\emptyset$ and $N_2(y_0)\cap R_0\neq\emptyset$. This defines an alternating cycle $x_0,u_0,P_1,b_0,a_0,P_2,v_0,y_0,x_0$. Here $u_0$ is a $c_2$-neighbor of $x_0$ in $\hT_\ell$ and $P_1$ is  (the reversal of) a path from $u_0$ to $b_0$ and $P_2$ is the path from $a_0$ to $v_0\in X_k$, $v_0\in N_2(y_0)$. This completes the proof of Theorem \ref{th1}.

{\bf Verification of \eqref{cl1}, \eqref{cl2}}:  We have by the assumption $a_0\in R_0$ that 
\[
|X_1|=|Y_1|\geq \frac{\a_2\b\log n}{10}-o(\log n)
\]
Now suppose that $1\leq |X_i|\leq n/(200\log n)$. Then, by \eqref{props},
\[
e_2(X_i:(N_2(X_i)\setminus M(A_2\setminus R_0)))\geq  \frac{(\a_2\b\log n)|X_i|}{10+o(1)}.
\]
Applying Lemma \ref{lem1}(a) we see that 
\beq{cl3}{
|N_2(X_i)\setminus M(A_2\setminus R_0)|\geq\frac{(\a_2\b\log n)|X_i|}{20+o(1)}.
}
Because the sets $X_1,X_2,\ldots$ expand rapidly, the total size of $\bigcup_{j\leq i}X_j$ is small compared with the R.H.S of \eqref{cl3} and \eqref{cl1} follows. The argument for \eqref{cl2} is similar.
\proofend
\section{Concluding Remarks}
We have established that w.h.p. $mcp(G)$ is almost all of $[0,n]^q$ and posed the question of findng the exact threshold for $mcp(G)=[0,n]^q$. It seems technically feasible to extend our results to randomly colored $G_{n,p}$. We leave this for future research. It would be of some interest to analyse other spanning subgraphs from this point of view e.g. Hamilton cycles.

\end{document}